\documentclass[12pt]{article}
\usepackage[latin1]{inputenc}
\usepackage{amssymb,amsmath,amsfonts}
\newtheorem{theorem}{Theorem}

\newtheorem{lemma}[theorem]{Lemma}

\newcommand{\R}{\mathbb{R}}

\newcommand{\Sf}{\mathbb{S}}
\newcommand{\C}{\mathbb{C}}

\newcommand{\spa}{\mbox{span}}

\newcommand{\po}{{\hspace*{-1ex}}{\bf .  }}

\newcommand{\nap}{\nabla^{\perp}}

\newcommand{\nab}{\tilde\nabla}

\newcommand{\la}{\lambda}
\def\<{{\langle}}
\def\>{{\rangle}}

\def\J{{\cal J}}
\def\n{\nabla}
\def\d{\partial}

\def\a{\alpha}

\def\bea{\begin{eqnarray*} }
\def\eea{\end{eqnarray*} }
\def\be{\begin{equation} }
\def\ee{\end{equation} }

\def\nap{\nabla^\perp}
\def\proof{\noindent\emph{ Proof: }}
\def\qed{\ifhmode\unskip\nobreak\fi\ifmmode\ifinner
\else\hskip5 pt \fi\fi\hbox{\hskip5 pt \vrule width4 pt
height6 pt  depth1.5 pt \hskip 1pt }}

\begin{document}

\title{A class of superconformal surfaces}
\author{M. Dajczer and Th. Vlachos}
\date{}
\maketitle

\begin{abstract} Superconformal surfaces in Euclidean space are the ones for
which the ellipse of curvature at any point is a nondegenerate circle. 
They can be characterized as the surfaces for which a well-known pointwise inequality
relating the intrinsic Gauss curvature with the extrinsic normal and mean curvatures, 
due to Wintgen (\cite{Wi}) and Guadalupe-Rodr\'iguez (\cite{GR}) for any codimension,
reaches equality at all points. 
In this paper, we show that any pedal surface to a $2$-isotropic Euclidean surface is 
superconformal. Opposed to almost all known examples, superconformal surfaces in this class 
are not conformally equivalent to minimal surfaces. Moreover, they can be given in an 
explicit parametric form since $2$-isotropic surfaces admit a Weierstrass type 
representation.
\end{abstract}

A  surface in Euclidean space  $f\colon M^2\to\R^n$, $n\geq 4$, 
is called \emph{superconformal}  if at any point the ellipse 
of curvature is a nondegenerate circle.  The \emph{ellipse of 
curvature} ${\cal E}(p)$ of $f$ at $p\in M^2$ is the ellipse in the normal 
space $N_fM(p)$ given as the image of the second fundamental form 
$\a_f$ restricted to the unit circle in the tangent plane $T_pM$, i.e., 
$$
{\cal E}(p)=\{\a_f(X,X):X\in T_pM,\,\|X\|=1\};
$$
see \cite{GR} and references 
therein for several facts on this concept whose study started almost a century 
ago due to the work of Moore and Wilson \cite{mw1}, \cite{mw2}.  
Superconformality is invariant under conformal (Moebius) transformations since 
the property of ${\cal E}(p)$ being a circle is invariant under conformal changes 
of the metric  of the ambient space. 

Perhaps the most notable property of superconformality is that the surfaces 
yield equality in the pointwise inequality due to  Wintgen \cite{Wi} for
$n=4$ and Guadalupe- Rodr\'iguez \cite{GR} for any codimension, that reads 
as follows:
$$
K + |K_N|\leq \|H\|^2
$$
where $K$ denotes the Gaussian curvature of $M^2$ and $K_N$ and $H$ are the 
normal curvature and mean curvature vector of $f$, respectively. 

In the special case of surfaces $f\colon M^2\to\R^4$ the condition of superconformality 
is rather strong. This allowed Dajczer and Tojeiro \cite{DT} to provide a 
parametric representation of all non-minimal superconformal surface in terms 
of any minimal surface and its conjugate. Notice that minimal superconformal surfaces 
in $\R^4$ are just holomorphic curves with respect to some complex structure in 
the ambient space. 

In higher codimension superconformality is not longer such a strong assumption
and several classes of examples of such surfaces have been considered.
The ``trivial" examples are the ones conformally equivalent (i.e., congruent 
by a Moebius transformation) to  minimal superconformal surfaces. These 
Euclidean surfaces are called  \emph{$1$-isotropic} and admit a Weierstrass type 
representation as, in fact, does any simply connected \mbox{$s$-isotropic} surface; 
see next section for definitions and some properties including a parametric 
representation of these minimal surfaces.

Other examples of superconformal surfaces are the images by stereographic 
projection of minimal superconformal surfaces in the sphere and hyperbolic space.
In the spherical case, this class of surfaces has been studied in different
contexts; see \cite{BPW}, \cite{Mi} and \cite{Vl}. In particular, it 
includes the minimal $2$-spheres studied in \cite{Ba}, \cite{Cal} and  \cite{Ch}, the
holomorphic curves in the nearly Kaehler sphere $\Sf^6$ \cite{Br} and Lawson's surfaces
\cite{Vl}  for appropriate choice of the several parameters involved in their 
definition. There are also the super Willmore surfaces in even 
dimensional spheres introduced as isotropic \mbox{S-Willmore} surfaces by Ejiri \cite{Ej} 
and classified in terms of isotropic holomorphic curves in complex projective spaces. 
It was shown recently \cite{DV} that superconformal S-Willmore surfaces are the ones to have 
a dual  surface in the same class. Finally, the case of the hyperbolic ambient space was 
studied in \cite{Hu}.

All of the above ``nontrivial" examples have in common that they do not posses an
explicit parametric description.  Moreover, with the exception of the super Willmore 
surfaces, they are all conformally  equivalent to minimal surfaces.  
In this paper, we present a new class of superconformal surfaces which do admit 
an explicitly  parametric representation and are never locally conformally equivalent 
to minimal surfaces.\medskip

The \emph{pedal surface} $g\colon M^2\to\R^n$  to a surface 
$f\colon M^2\to\R^n$ is defined as the locus of the foot of the perpendicular 
from a fixed arbitrary point to the tangent planes to $f$ and  assumed to be 
an immersion.   For surfaces in $\R^3$ this concept was considered by several authors 
during the second  half of  XIX century; see  \cite[p.\ 673]{en} and 
\cite[pp.\ 155-161]{so}.

\begin{theorem}\po\label{main}  
A pedal surface $g\colon M^2\to\R^n$, $n\geq 6$, to a  substantial $1$-isotropic surface 
$f\colon M^2\to\R^n$ is superconformal if and only $f$ is $2$-isotropic. 
Then, we have:
\begin{itemize}
\item[(i)] The surfaces $f$ and $g$ are conformal.
\item [(ii)] The surface $g$ is not locally conformally equivalent to 
a minimal surface.
\item [(iii)] The surface $g$ is never S-Willmore.
\end{itemize}
\end{theorem}

In particular, the above result applies to any holomorphic curve $f\colon M^2\to\C^m$ for 
$m\geq 3$,  since these surfaces are the ones that are $(m-1)$-isotropic in $\R^{2m}$. By
$f$ being substantial we mean that the surface is not contained in any proper affine
subspace.
\vspace{1ex}

It was shown by Rouxel \cite{Ro} that a pedal surface to any holomorphic curve 
in $\C^2\cong\R^4$ is a superconformal surface. But opposed to part $(ii)$ in the theorem he 
also proved that the pedal surface is conformally equivalent to an holomorphic curve 
with respect to some conformal structure in $\R^4$.  Besides, the proof of these two facts 
can be easily obtained from the arguments in this paper.

\section{Preliminaries}

Let $f\colon M^2\to\R^n$, $n\geq 4$, be an isometric immersion of a
two-dimensional Riemannian manifold into $n$-dimensional Euclidean space. 
The $k^{th}$\emph{-normal space} of the surface $f$ at $p\in M^2$ for $k\geq 1$ 
is defined as
$$
N^f_k(p)=\spa\{\alpha_f^{k+1}(X_1,\ldots,X_{k+1}):X_1,\ldots,X_{k+1}\in T_pM\}
$$
where $\alpha_f^2=\alpha_f\colon TM\times TM\to N_fM$ is the standard second 
fundamental form and  
$$
\alpha_f^s\colon TM\times\cdots\times TM\to N_fM,\;\; s\geq 3, 
$$
denotes the symmetric tensor called
the $s^{th}$\emph{-fundamental form} defined inductively by
$$
\alpha_f^s(X_1,\ldots,X_s)=\left(\nabla^\perp_{X_s}\ldots
\nabla^\perp_{X_3}\alpha_f(X_2,X_1)\right)^\perp
$$
where $(\;\;)^\perp$ means taking the projection onto
the normal subspace $(N^f_1\oplus\ldots\oplus N^f_{s-2})^\perp$.

Since this paper deals with local theory, we always assume that the immersion 
$f$ is \emph{regular} (called nicely curved in \cite{Sp}) which means  that all 
$N^f_k$'s have constant dimension for each $k$ and thus form normal subbundles. 
Notice that this condition is verified along connected components of an open 
dense subset of $M^2$.

Assume now that the surface $f$ is minimal. In this case, the normal bundle 
splits as
\be\label{decon}
N_fM=N_1^f\oplus N_2^f\oplus\dots\oplus N_m^f,\;\;\; m=[(n-1)/2],
\ee
since all higher normal bundles  have rank two except the last one that has 
rank one if $n$ is odd; see \cite{Ch}, \cite{DF} or \cite{Sp}. 
Moreover, the orientation of $M^2$ and $\alpha_f^s$ induce an orientation 
on each plane bundle $N_{s-1}^f$ given by the ordered pair
$$
\{\a_f^s(X,\ldots,X),\;\a_f^s(JX,\ldots,X)\}
$$
where $J$ is the complex structure of $M^2$ determined by the metric and 
the orientation. 

Each two dimensional normal subbundle in (\ref{decon}) has associate an ellipse 
of curvature defined  as follows: The \textit{$s^{th}$-order ellipse of curvature}
$\mathcal{E}_s^f(p)\subset N^f_s(p),\, 1\leq s\leq m_0$,\, of $f$ at $p\in M^2$ is
$$
\mathcal{E}_s^f(p)=\{\alpha_f^{s+1}(X_{\psi},\dots,X_{\psi}) :
X_{\psi}=\cos\psi X+\sin\psi JX\;\;\mbox{and}\;\;\psi\in [0,2\pi)\}
$$
where $X\in T_p M$ has  unit length and $N_{m_0}^f$ is the last plane subbundle. 
Of course, the standard ellipse of curvature is just $\mathcal{E}(p)=\mathcal{E}_1^f(p)$.
\vspace{1ex}

The minimal surface $f$ is called \emph{$m$-isotropic} if
$\mathcal{E}_r^f(p)$ is a circle for any $p\in M^2$ and any $1\leq r\leq m$.
Thus minimal superconformal surfaces are the $1$-isotropic ones. A fundamental fact in
relation to this paper is that any simply connected $m$-isotropic surface admit a 
Weierstrass type representation  given in \cite{DG} based on results in \cite{CCh}.

Start with a nonzero
holomorphic map 
$$
\alpha_0\colon\, U\to\C^{N-2(m+1)}.
$$ 
Assuming that 
$$
\alpha_r\colon\, U\to\C^{N-2(m-r+1)},\;\;0\leq r\leq m,
$$ 
has been defined already, set
$$
\alpha_{r+1}=\beta_{r+1}\left(1-\phi_r^2,i(1+\phi_r^2),2\phi_r\right)
$$
where $\phi_r=\int^z\alpha_rdz$ and $\beta_{r+1}\neq 0$ is any holomorphic 
function. Then, the surface $g=\mbox{Re}\,\phi_{m+1}$ 
in $\R^N$ is  $m$-isotropic.

\section{The proof}

Given a surface $f\colon M^2\to\R^n$, $n\geq 4$, set $\a_{ij}=\a_f(X_i,X_j)$
where $\{X_1,X_2\}$ is an orthonormal basis  of $T_pM$ at $p\in M^2$.
For a unit vector $v=\cos\theta X_1+\sin\theta X_2$, we have 
$$
\a_f(v,v)=H+\cos2\theta\,\xi_1+\sin2\theta\,\xi_2
$$
where $\xi_1=\frac{1}{2}(\a_{11}-\a_{22})$, $\xi_2=\a_{12}$
and $H$ is the mean curvature vector. 
Thus ${\cal E}(p)$ is a circle if and only if we have
\be\label{super}
\<\a_{12},\a_{11}-\a_{22}\>=0
\mbox\;\;\;\mbox{and}\;\;\;\|\a_{11}-\a_{22}\|=2\|\a_{12}\|.
\ee

The complexified tangent bundle $TM\otimes \mathbb{C}$ decomposes into
the eigenspaces of the complex structure $J$ corresponding to the eigenvalues 
$i$ and $-i$ denoted by $T^{\prime}M$ and $T^{\prime\prime}M$, respectively. 
 The second fundamental form   can
be complex linearly extended to $TM\otimes \mathbb{C}$ with values in the
complexified vector bundle $N_f M \otimes \mathbb{C}$ and then decomposed
into its $(p,q)$-components, $p+q=2,$ which are tensor products of $p$ many
1-forms vanishing on $T^{\prime \prime }M$ and $q$ many 1-forms vanishing 
on $T^{\prime }M$.

Taking local isothermal coordinates $\{x,y\}$  and $\d=(1/2)(\d_x-i\d_y)$,
we have from (\ref{super}) that $f$ is superconformal if and only if the 
$(2,0)$-part of the second fundamental form is isotropic, i.e.,  the 
complex line bundle spanned by $\a_f(\d,\d)$ is isotropic.
Moreover, the surface is minimal if its $(1,1)$-part of $\a_f$
vanishes, i.e., if $\a_f(\d,\bar\d)=0$.

Without loss of generality, we choose to take the pedal surface always with 
respect to the origin since changing that amounts to add a constant vector 
to $f$; see last section for further discussion on this election. 

Given a substantial minimal surface $f\colon M^2\to\R^n,\; n\geq 6$,  
decompose the position vector into its tangent and normal components
\be\label{deco}
f=f_*Z+g
\ee
and assume that $g\colon M^2\to\R^n$ is an immersion. 
Thus $g$ is the pedal surface of $f$ with respect to the origin.
Differentiating (\ref{deco}) yields
\be\label{d}
X=\n_X Z-A_g X\;\;\mbox{and}\;\;\a(X,Z)+\nap_Xg=0
\ee\label{Z}
where  $A_g=A^f_g$ is the shape operator associated to $g$ and $\a=\a_f$. 
Moreover,  we had identified $X$ with $f_*X$ for $X\in TM$. It follows that
\be\label{g*}
g_*X=-A_g X-\a(X,Z).
\ee

Hereafter, we assume that $Z\neq 0$ and $g\not\in (N_1^f)^\perp$ at any point. 
Clearly, we cannot have $Z=0$ in an open subset. On an open 
subset where $g\in (N_1^f)^\perp$ it is easy to see from 
the first equation in (\ref{d}) that $M^2$  would be flat and $f$ totally geodesic.
We denote $g=\delta+\eta$
where $\delta\neq 0$ and $\eta$ are the $N_1^f$-component and ${(N_1^f)^\perp}$-component 
of $g$, respectively.

\begin{lemma}\po\label{conf}  
A minimal surface  $f\colon M^2\to\R^n$ is conformal to its  pedal surface
if and only if it is $1$-isotropic. In this case, 
we have $ds_g^2=-(1/2)K\theta ds_f^2$
where $K$ is the Gaussian curvature of $M^2$ and 
$\theta=\|Z\|^2+\|\delta\|^2$.
\end{lemma}

\proof In general, we have that $g$ is conformal to $f$ if and only if $\<g_*\d,g_*\d\>=0$. 
Since $f$ is minimal, we obtain from (\ref{g*}) that
\bea
\<g_*X,g_*Y\>=-\det A_\delta\, \<X,Y\>+\<\a(X,Z),\a(Y,Z)\>
\eea
for any $X,Y\in TM$. Thus $g$ is conformal to $f$ if and only if
$\<\a(\d,Z),\a(\d,Z)\>=0$. Being $f$  minimal, we have 
$i\a(\d,Z)=\a(\d,JZ)$. Hence $g$ is conformal to $f$ if and only if
$$
\<\a(\d,Z-iJZ),\a(\d,Z-iJZ)\>=0,
$$
and this is equivalent to $f$ being $1$-isotropic since 
$Z-iJZ \in T'M={\mathrm{span}}_{\C}\{\d\}$.

We compute the conformal factor. We have
$$
\<\d,\bar\d \>_g =-\det A_\delta\, \<\d,\bar\d\>+\<\a(\d,Z),\a(\bar\d,Z)\>.
$$
 From the Gauss equation and since $f$ is minimal, it follows that
$$
\<\a(\d,Z),\a(\bar\d,Z)\>=-\<R(\d,Z)Z,\bar\d\>
=K(|\<\d,Z\>|^2-\|Z\|^2\<\d,\bar\d\>)=-\frac{1}{2}K\|Z\|^2\<\d,\bar\d\>.
$$
where $R$ is the curvature tensor of $M^2$. 
Since $f$ is $1$-isotropic, given an orthonormal $\{e_3,e_4\}$ frame in $N_1^f$ 
there is an orthonormal tangent frame such that 
\be\label{form}
A_{e_3}=\begin{pmatrix} 
\kappa&0\\ 
0&-\kappa&
\!\!\!\!\!\end{pmatrix},\;\;
A_{e_4}=
\begin{pmatrix} 
0& \kappa\\ 
\kappa& 0
\end{pmatrix}.
\ee
Thus $\det A_\delta=(1/2)K\|\delta\|^2$, and the claim follows.\qed

\begin{lemma}\po\label{Ng} Assume that  $f$ is $1$-isotropic. 
Then, the normal bundle of $g$ is given by
$$
N_gM=\spa\,\{Z-\delta,JZ+J_1^\perp\delta\}\oplus (N_1^f)^\perp
$$ 
where  $J_1^\perp$ denotes the orthogonal 
complex structure of $N_1^f$. 
\end{lemma}

\proof That $f$ is $1$-isotropic is equivalent to
$$
J_1^\perp\a(X,Y)=\a(JX,Y)
$$
for any $X,Y\in TM$. Then
$$
\<g_*X,JZ+J_1^\perp\delta\>=\<JA_\delta X+J_1^\perp\a(X,Z),Z+\delta\>
=-\<A_\delta JX,Z\>+\<\a(JX,Z),\delta\>=0,
$$
and the remaining cases are immediate.\qed
\vspace{1,5ex}

We define a complex structure
$\J= J_1^\perp\oplus J_2^\perp$ on $N_1^f\oplus N_2^f$ 
where  $J_2^\perp$ denotes the orthogonal 
complex structure of  $N_2^f$. 

\begin{lemma}\po\label{J} If $f$ is $1$-isotropic, then the following equations hold:
$$
(\nap_X \J)e_3= (\la-1)(*\omega(X)e_5+\omega(X)e_6),\;\;
(\nap_X \J)e_4= (\la-1)(\omega(X)e_5-*\omega(X)e_6)
$$
where $0<\la\leq 1$ is the ratio of the length of the axes 
of $\mathcal{E}_2^f$. Moreover,
the orthonormal frames $\{e_3, e_4\}$ and $\{e_5, e_6\}$ span $N_1^f$ and $N_2^f$, 
respectively, $\omega=\<\nap e_3,e_5\>$  and $*$ is the Hodge operator.
\end{lemma}

\proof Arguing as in the proof of Lemma~5 in \cite{Vl} we obtain that
there is a local orthonormal frame $\{ e_A\}, 1\leq A\leq n$, 
where $\{e_1,e_2\}$ is a tangent frame and $\{e_{2r+1},e_{2r+2}\}$ is a positively 
oriented frame of $N_r^f$ for $1\leq r \leq m_0$, such that the normal
connection forms $\omega_{\a\beta}=\< \nap e_{\a},e_{\beta}\>$
satisfy
\be\label{conn}
\omega_{45}=-*\omega,\;\;\omega_{46}=-*\omega_{36},\;\;\omega_{36}
=\la*\omega\;\;\mbox{and}\;\;\omega_{46}=\la*\omega_{45}.
\ee
Moreover, we have that $(\nap_X \J)e_{\a}\in N_2^f$ for $3\leq\a\leq 4$, 
and the remaining of the proof is straightforward.\qed

\begin{lemma}\po\label{ag} If  $f$ is $1$-isotropic, 
then the second fundamental form of $g$ satisfies
\bea
\a_g(\d,\d)\!\!\!&=&\!\!\!\<\a_g(\d,\d),\xi_1\>(\xi_1+i\xi_2)
-\omega(\d)\< \a(\d,Z), e_3+ie_4\> (e_5-i\lambda e_6)
\eea
where $\xi_1= (Z-\delta)/\sqrt{\theta}$\, and\, $\xi_2
=(JZ+\J\delta)/\sqrt{\theta}.$
\end{lemma}

\proof From
$$
\nab_Xg_*Y=-\n_XA_\delta Y+A_{\a(Y,Z)}X-\a(X,A_\delta Y)-\nap_X\a(Y,Z)
$$
we obtain that
\be\label{normals}
(\a_g(X,Y))^{(N_1^f)^\perp}=-(\nap_X\a(Y,Z))^{(N_1^f)^\perp}.
\ee
Hence, $\<\a_g(\d,\d),\nu\>=-\<\nap_\d\a(\d,Z),\nu\>$
for any $\nu\in (N_1^f)^\perp$. In particular, 
$$
\<\a_g(\d,\d),\nu\>=0\;\;\mbox{for any}\;\;\nu\in (N_1^f\oplus N_2^f)^\perp.
$$
On the other hand, since $\a(\d,A_\delta\d)=0$, we have
\be\label{agxi}
\<\a_g(\d,\d),Z-\delta\>=\<-\n_{\d}A_\delta\d+A_{\a(\d,Z)}\d,Z\>
+\<\nap_{\d}\a(\d,Z),\delta\>
\ee
and, using Lemma \ref{J}, that
\bea
\<\a_g(\d,\d),JZ+\J\delta)\>\!\!\!&=&\!\!\!
-\<\n_{\d}A_\delta\d-A_{\a(\d,Z)}\d,JZ\>
-\<\nap_{\d}\a(\d,Z),\J\delta\>\\
\!\!\!&=&\!\!\!\<J\n_{\d}A_\delta\d,Z\>+\<A_{\a(\d,Z)}J\d,Z\>
+\<\J\nap_{\d}\a(\d,Z),\delta\>\\
\!\!\!&=&\!\!\!-\<\n_{\d}A_\delta J\d,Z\>+i\<A_{\a(\d,Z)}\d,Z\>
+\<\nap_{\d}\J\a(\d,Z), \delta \>\\
\!\!\!&=&\!\!\!-i\<\n_{\d}A_\delta\d,Z\>+i\<A_{\a(\d,Z)}\d,Z\>
+i\<\nap_{\d}\a(\d,Z),\delta\>\\
\!\!\!&=&\!\!\!i\<\a_g(\d,\d),Z-\delta\>.
\eea
It follows that
$$
\a_g(\d,\d)= \<a_g(\d,\d),\xi_1\>(\xi_1+i\xi_2)
+\<\a_g(\d,\d),e_5\>e_5  + \<\a_g(\d,\d),e_6\>e_6.
$$
Moreover, we obtain using (\ref{conn}) that
\begin{eqnarray}\label{rel}
\<\a_g(\d,\d),e_5\>=\<\a(\d,Z),\nap_\d e_5\>
\!\!\!&=&\!\!\!-\<\a(\d,Z),\omega(\d)e_3-*\omega(\d)e_4\>\nonumber\\
\!\!\!&=&\!\!\!-\omega(\d)\<\a(\d,Z),e_3+i e_4\>.
\end{eqnarray}
Similarly, 
$$
\<\a_g(\d,\d),e_6\>=i\lambda\omega(\d)\<\a(\d,Z) e_3+ie_4\>,
$$
and the result follows.\qed
\vspace{1,5ex}

\noindent\emph{Proof of Theorem \ref{main}:} By Lemma \ref{ag} 
the pedal surface is superconformal if and only if
$$
(1-\lambda^2)(\omega(\d) \<\a(\d,Z),e_3+ie_4\>)^2=0.
$$
But if $\<\a(\d,Z),e_3+ie_4\>=0$, then  $Z=0$ and that is a contradiction.
And if $\omega(\d)=0$ we have from (\ref{conn}) that $N_1^f$ is parallel 
in the normal bundle, also a contradiction. We conclude that $\lambda=1$.

 From  Lemma \ref{conf}, the mean curvature 
vector field of $g$ is 
$$
H_g=-\frac{1}{K\theta}\Delta g
$$
where the Laplace operator of $M^2$  acts on vector valued functions. 
Using (\ref{g*}), we  have
$$
\Delta g=-\sum_{j=1,2}((\n_{e_j} A_{\delta})e_j 
- A_{\a(e_j,Z)}e_j+\a(e_j, A_{\delta}e_j)+\nap_{e_j}\a(e_j,Z)-\a(\n_{e_j}e_j,Z))
$$
for an orthonormal frame $\{e_1,e_2\}$ in $TM$. 
Using the Codazzi equation, we have
\bea
\sum_{j=1,2}\<(\n_{e_j} A_{\delta})e_j,X\>\!\!\!&=&\!\!\!
\sum_{j=1,2}\<(\n_{e_j} A_{\delta})X,e_j\> \\
\!\!\!&=&\!\!\!\sum_{j=1,2}\<(\n_{X} A_{\delta})e_j,e_j\>+ \sum_{j=1,2}
\<A_{\nap_{e_j}\delta}X,e_j\>
-\sum_{j=1,2} \<A_{\nap_{X}\delta}e_j,e_j\>\\
\!\!\!&=&\!\!\!\sum_{j=1,2}\<(\n_{X} A_{\delta})e_j,e_j\>- \sum_{j=1,2}
\<A_{\a(e_j,Z)}X,e_j\>.
\eea
On the other hand,
\bea
\sum_{j=1,2}\<(\n_{X} A_{\delta})e_j,e_j\>\!\!\!&=&\!\!\!\sum_{j=1,2}\<\n_{X}
A_{\delta}e_j,e_j\>
-\sum_{j=1,2}\< A_{\delta}\n_{X}e_j,e_j\>\\
\!\!\!&=&\!\!\!\sum_{j=1,2}X\< A_{\delta}e_j,e_j\>
-2\sum_{j=1,2}\<A_{\delta}\n_{X}e_j,e_j\>\\
\!\!\!&=&\!\!\!0.
\eea
We obtain that
\be\label{delta}
\Delta g=\sum_{j=1,2}(2 A_{\a(e_j,Z)}e_j-\a(e_j,A_{\delta}e_j)  
- \nap_{e_j}\a(e_j,Z) +\a(\n_{e_j}e_j, Z)).
\ee
Choosing orthonormal frames $\{e_3,e_4\}$  in $N_1^f$ with $e_3=\delta/\|\delta\|$ 
and $\{e_1,e_2\}$ in $TM$ such that (\ref{form}) holds, we have
$$
\sum_{j=1,2}A_{\a(e_j,Z)}e_j=-KZ\;\;\text{and}\;\; 
\sum_{j=1,2}\a(e_j, A_{\delta}e_j) = -K\delta.
$$
A straightforward computation using (\ref{d}), (\ref{conn}) and 
the Codazzi equation for $A_{e_3}$, that is,
$$
e_1(\kappa)= -2\kappa\psi(e_2)+\kappa\omega_{34}(e_2)\;\;\mbox{and}\;\;
e_2(\kappa)= 2\kappa\psi(e_1)-\kappa\omega_{34}(e_1),
$$
gives that
$$
\sum_{j=1,2}\nap_{e_j}\a(e_j,Z)=-K\delta-\kappa((z_1\psi(e_2)+z_2\psi(e_1))e_3
-(z_1\psi(e_1)-z_2\psi(e_2))e_4)
$$
where $Z=z_1e_1+z_2e_2$ and $\psi=\<\n e_1,e_2\>$.  Moreover, 
$$
\sum_{j=1,2}\a(\n_{e_j}e_j,Z)=-\kappa((z_1\psi(e_2)+z_2\psi(e_1))e_3
-(z_1\psi(e_1)-z_2\psi(e_2))e_4).
$$
Now, we obtain from (\ref{delta}) that
$\Delta g=2K(\delta-Z)$ and, consequently, that 
\be\label{meang}
H_g= \frac{2}{\theta}(Z-\delta).
\ee

Consider the inversion $\mathcal{I}$  with respect to a sphere with 
radius $R$ centered at $p_0$ and the immersion 
$\tilde{g}=\mathcal{I}\circ g$.  Then, there is  a vector bundle isometry
$\mathcal{P}$ between the normal bundles $N_gM$ and $N_{\tilde{g}}M$
(see \cite{DT}) given by
$$
\mathcal{P}\mu=\mu-2\frac{\<g-p_0,\mu\>}{\|g-p_0\|^2}(g-p_0)
$$
such that shape operators of $g$ and $\tilde{g}$ are related by
$$
\tilde{A}_{\mathcal{P}\mu}=\frac{1}{R^2}\left(\|g-p_0\|^2A_\mu+2\<g-p_0,\mu\>I \right).
$$
Then, the mean curvature vector of $\tilde{g}$ is given by
$$
H_{\tilde{g}}=\frac{1}{R^2}\mathcal{P}(\|g-p_0\|^2 H_g+2(g-p_0)^\perp).
$$
Using Lemma \ref{Ng}, we have
\bea
\frac{1}{2}\|g-p_0\|^2 H_g+(g-p_0)^\perp \!\!\!&=&\!\!\!\frac{1}{\theta}(\|g-p_0\|^2
-\|\delta\|^2-\<p_0 ,Z-\delta\>)(Z-\delta)\\
\!\!\!&&\!\!\!-\frac{1}{\theta}\<p_0,JZ+\J\delta\>(JZ+\J\delta)
+ \eta-(p_0)^{(N_1^f)^\perp}.
\eea
Thus $\tilde{g}$ is minimal if and only if
\be\label{min}
\|g-p_0\|^2-\|\delta\|^2-\<p_0 ,Z-\delta\>=0,\; \<p_0 ,JZ+\J\delta\>=0 \;\; 
\text{and}\;\; \eta=(p_0)^{(N_1^f)^\perp}. 
\ee
If $p_0=0$, then $g \in N_1^f$ and (\ref{d}) 
combined with Lemma \ref{J} yield that $N_1^f=\spa\{\delta,\J\delta\}$ 
is parallel in the normal bundle,  which has been excluded.
Thus, we may assume $p_0\neq 0$. Setting $p_0= W+ \xi+\eta$ where $W \in TM$ 
and $\xi \in N_1^f$, we have
\be\label{w}
\n_XW=A_{\xi}X\;\;\mbox{and}\;\;\a(X,W)=-\nap_X\xi-\nap_X\eta.
\ee
 From (\ref{d}) we obtain $\nap_X (\delta-\xi) \in N_1^f$ and if $\delta\neq\xi$
the same argument as before gives that $N_1^f$ is parallel in the normal bundle. 
Hence  $\xi=\delta$. Hence, the second equation in (\ref{min}) gives that $W$ and $Z$ 
are collinear. But then the first equation in (\ref{min})  yields  $W=Z$, 
and this is not possible by (\ref{d}) and (\ref{w}). Thus $\tilde g$ cannot be 
minimal.\vspace{1ex} 

To conclude the proof we show that $g$ is never S-Willmore. 
By definition, we have that $g$ is S-Willmore if and only if 
$\hat{\nabla}_{\bar{\d}}^\perp\a_g(\d,\d)$ is parallel to $\a_g(\d,\d)$ 
where $\hat{\nabla}^\perp $ is the normal connection of $g$. 
The Codazzi equation gives
$$
\hat{\nabla}_{\d}^\perp H_g 
=\frac{2}{ \rho^2}\hat{\nabla}_{\bar{\d}}^\perp\a_g(\d,\d)
$$
where the induced metric is $ds_g^2=\rho^2 |dz|^2$ and $z=x+iy$.
Therefore, the superconformal surface $g$ is S-Willmore if and only if 
$\hat{\nabla}_{\d}^\perp H_g$ is parallel to $\a_g(\d,\d)$.

By Lemma \ref{ag}, we have  
$$
\a_g(\d,\d)=\<\a_g(\d,\d),\xi_1\>(\xi_1+i\xi_2)
-\omega(\d)\< \a(\d,Z), e_3+ie_4\> (e_5-i\lambda e_6).
$$
We claim that
$$
\<\a_g(\d,\d),\xi_1\>=\frac{1}{\sqrt{\theta}}(-\d(\|\delta\|)\<A_{e_3}\d,Z\> 
+\|\delta\|\<A_{e_3}\d,\d\>-\|\delta\|\omega_{34}(\d)\<A_{e_4}\d,Z\>).
$$
First observe that since $f$ is 1-isotropic, we have
$$
\<A_{\a(\d,Z)}\d ,Z\>=\< A_{e_3}\d,Z\>^2+ \< A_{e_4}\d,Z\>^2=0.
$$
Now from (\ref{agxi}) using (\ref{d}), we obtain
\bea
\<\a_g(\d,\d),\xi_1\>
\!\!\!&=&\!\!\!\frac{1}{\sqrt{\theta}}(-\<\n_\d A_\delta \d,Z\>
+\<\nap_\d \a(\d,Z),\delta\>)\\
\!\!\!&=&\!\!\!\frac{1}{\sqrt{\theta}}(-\<\n_\d A_\delta\d,Z\>
+\d\<\a(\d,Z),\delta\>-\<\a(\d,Z),\nap_\d\delta\>)\\
\!\!\!&=&\!\!\!\frac{1}{\sqrt{\theta}}(-\<\n_\d A_\delta \d,Z\>
+\d \<A_\delta \d,Z\>-\<\a(\d,Z),\nap_\d (\|\delta\|e_3)\>)\\
\!\!\!&=&\!\!\!\frac{1}{\sqrt{\theta}}(\<A_\delta\d,\d\>
+\<A_\delta\d,A_\delta\d\>-\d(\|\delta\|)\<A_{e_3}\d,Z\>
-\|\delta\|\omega_{34}(\d)\<A_{e_4}\d,Z\>),
\eea
and the claim follows.

On the other hand, the second equation in (\ref{d}) using (\ref{conn}) yields
$$
\<A_{e_3}\d,Z\>=-\d \|\delta\|+\omega(\d)(\eta_5-i\eta_6)\;\; \text{and} \;\; 
\<A_{e_4}\d,Z\>=- \|\delta\|\omega_{34}(\d)+i\omega(\d)(\eta_5-i\eta_6)
$$
where $\eta_j=\<\eta,e_j\>,\,j=5,6$.
We obtain that
\begin{eqnarray}\label{one}
\a_g(\d,\d)\!\!\!&=&\!\!\!\frac{1}{\sqrt{\theta}}(\|\delta\|\<A_{e_3}\d,\d\>
-\omega(\d)\< \a(\d,Z),e_3+ie_4\>(\eta_5-i\eta_6))(\xi_1+i\xi_2)\nonumber\\
\!\!\!&&\!\!\!-\omega(\d)\< \a(\d,Z),e_3+ie_4\>(e_5-ie_6).
\end{eqnarray}
On account of (\ref{meang}) and Lemma \ref{Ng}, we obtain
$$
\hat{\nabla}_{\d}^\perp H_g= -\frac{2}{\theta^{3/2}}\d\theta\xi_1 
+ \frac{2}{\theta}(\nab_\d (Z-\delta))^{N_gM}.
$$
Using (\ref{d}) we find that
$$
\nab_\d (Z-\delta)=\d+2A_\delta \d +2\a(\d,Z)+\nap_\d \eta.
$$
It follows that
$$
\<\nab_\d (Z-\delta), \xi_1 \>
=\frac{1}{\sqrt \theta}(\<\d,Z\>+\<\nap_\d \delta, \eta\>)
$$
and 
$$
\< \nab_\d (Z-\delta), \xi_2 \>
=-\frac{1}{\sqrt \theta}(i\<\d,Z\>+\<\nap_\d \J\delta, \eta\>).
$$
Using (\ref{d}) we obtain that
$$
\d\theta=2\<\d,Z\>+2\<\nap_\d \delta, \eta\>
$$
and from (\ref{conn}) that
$$
(\nap_\d \eta)^{(N_1^f)^\perp}=(\nap_\d g)^{(N_1^f)^\perp}
-(\nap_\d \delta)^{(N_1^f)^\perp}
=-\|\delta\|\omega(\d)(e_5-ie_6).
$$
It follows easily  that
\begin{eqnarray}\label{two}
\hat{\nabla}_{\d}^\perp H_g\!\!\!&=&\!\!\!-\frac{2}{\theta^{3/2}} 
((\<\d,Z\>+\|\delta\|\omega(\d)(\eta_5-i\eta_6)))(\xi_1+i\xi_2)
-\frac{2}{\theta}\|\delta\|\omega(\d)(e_5-ie_6).
\end{eqnarray}
Hence, we conclude from (\ref{one}) and  (\ref{two}) that $g$ is S-Willmore if and only if
$$
\omega(\d)(\|\delta\|^2 \<A_{e_3}\d,\d\> +\<\d,Z\>\< \a(\d,Z), e_3+ie_4\>)=0.
$$
But $\omega(\d)\neq 0$ since, otherwise, we have from (\ref{conn}) that $N_1^f$ 
is parallel in the normal bundle, which has been  excluded. It follows easily 
that $g$ is S-Willmore if and only if $\kappa\,\theta=0$, but this is not possible.
\qed

\section{Final comments}

Changing the origin with respect to which the pedal of $f\colon M^2\to\R^n$  is
taken amounts to add a constant vector to $f$. Given $c\in\R$ and 
$v\in\R^n$, the pedal surface $g_{c,v}=(cf+v)^\perp$ of $cf+v$ is a superconformal 
surface which is conformal to $f$. Decomposing  $f$ and $v$ 
into its tangent and normal components
$$
f=Z+g\;\;\text{and}\;\;v=V+v^\perp,
$$
we have that $g_{c,v}=cg+v^\perp$. When $f$ is $2$-isotropic, it turns out that  
$g_{0,v}=v^\perp$ is superconformal which can be proved by a similar computations 
as the one given for $g$.  But in this case,  the image of $g_{0,v}=v^\perp$ by 
an inversion centered at $v$  is a minimal surface, hence
$1$-isotropic.
\vspace{2ex}

The superconformal surface $g$ constructed in Theorem \ref{main} 
not only satisfies part $(ii)$ but a much stronger condition, namely, 
for any surface  conformally equivalent to $g$ the rank of its first 
normal bundle is three. We sketch a proof of this fact but only in 
the special case of the pedal surface to a substantial $3$-isotropic surface 
$f$ in $\R^n$ with  $n\geq 8$.

Consider the immersion $\tilde{g}=\mathcal{I}\circ g$ where $\mathcal{I}$ is
the inversion with respect to a sphere of radius $R$ centered at $p_0$.
Then,
$$
\a_{\tilde{g}}(X,Y)=\frac{1}{R^2}\mathcal{P}(\|g-p_0\|^2 \a_g+2\<X,Y
\>(g-p_0)^\perp)
$$
where $\mathcal{P}\colon N_gM\to N_{\tilde{g}}M$ is the vector bundle isometry
already discussed. Then,
$$
N_1^{\tilde{g}}=\spa\{\mathcal{P}\mu_1,\mathcal{P}\mu_2,\mathcal{P}\mu_3\}
$$
where
\bea
\mu_1\!\!\!&=&\!\!\!\frac{1}{\sqrt{\theta}}((\|\delta\|
-\omega(Z)\eta_5-\omega(JZ)\eta_6)\xi_1
+(\omega(JZ)\eta_5-\omega(Z)\eta_6)\xi_2)-\omega(Z)e_5-\omega(JZ)e_6,\\
\mu_2\!\!\!&=&\!\!\!\frac{1}{\sqrt{\theta}}((\omega(Z)\eta_6-\omega(JZ)\eta_5)\xi_1
+(\|\delta\|-\omega(Z)\eta_5-\omega(JZ)\eta_6)\xi_2)-\omega(JZ)e_5+\omega(Z)e_6,\\
\mu_3\!\!\!&=&\!\!\!\frac{1}{\sqrt{\theta}} ((\|g-p_0\|^2-\|\delta\|^2
-\<p_0 ,Z-\delta\>)\xi_1-\<p_0 ,JZ+\J\delta\>\xi_2)+\eta-(p_0)^{(N_1^f)^\perp}.
\eea

Suppose that $\lambda_1\mu_1+\lambda_2\mu_2+\lambda_3\mu_3=0$
with $(\lambda_1,\lambda_2,\lambda_3)\neq (0,0,0)$. Then,
$$
\lambda_3 (\eta-p_0 )^{( N_1^f\oplus  N_2^f)^\perp}=0.
$$ 
If $\lambda_3=0$ then  $\omega=0$,  a contradiction. Thus $\lambda_3\neq 0$ 
and hence $\eta-p_0\in TM\oplus N_1^f\oplus  N_2^f$. We obtain that
$$
\a(X,W)=-\nap_X \xi-\nap_X\zeta
$$
where $p_0=W+\xi+\zeta,\, W\in TM,\,\xi\in N_1^f\,
\mbox{and}\,\zeta\in (N_1^f)^\perp$.

We claim that $\eta=\zeta$. Set $\eta=\zeta+\nu$ where $\nu\in N_2^f$.
The above combined with (\ref{d}) yield $\nap_X(\xi-\delta-\nu)\in N_1^f$
and, consequently, $\nap_X \nu \in N_1^f\oplus  N_2^f$. Since $f$ is 
$3$-isotropic, we have $\nap_X\J\nu=\J(\nap_X\nu)$. 
Thus, if $\nu\neq 0$, then  $N_1^f\oplus N_2^f$
is parallel in the normal bundle, a contradiction.
Hence $\eta=\zeta$,  $\nap_X (\xi-\delta)\in N_1^f$ and $\nap_X
\J(\xi-\delta)\in N_1^f$.
Since $f$ is substantial, we obtain that $\xi=\delta$ and, consequently, 
that $\a(X,W-Z)=0$. But this gives that $W=Z$, a contradiction.

\vspace{.5in} {\renewcommand{\baselinestretch}{1}
\hspace*{-20ex}\begin{tabbing} \indent\= IMPA -- Estrada Dona Castorina, 110
\indent\indent\= Univ. of Ioannina -- Math. Dept. \\
\> 22460-320 -- Rio de Janeiro -- Brazil  \>
45110 Ioannina -- Greece \\
\> E-mail: marcos@impa.br \> E-mail: tvlachos@uoi.gr
\end{tabbing}}

\begin{thebibliography}{lbl}


\bibitem{Ba} J. Barbosa, \textit{On minimal immersions of $S^2$
into $S^{2m}$.} Trans. Amer. Math. Soc. \textbf{210} (1975), 75--106.

\bibitem{BPW} J. Bolton, F. Pedit and L. Woodward, \emph{Minimal surfaces and
the affine Toda field model.} J. Reine Angew. Math. \textbf{459} (1995),
119--150.

\bibitem{Br} R. Bryant, \textit{Submanifolds and special structures on the
octonians.} J. Differential Geom. \textbf{17} (1982), 185--232.

\bibitem{Cal} E. Calabi, \emph{Minimal immersions of surfaces in Euclidean spheres.} 
J. Differential Geom. \textbf{1} (1967), 111--125.

\bibitem{CCh} C.\ C.\ Chen, \emph{The generalized curvature
ellipses and minimal surfaces.} Bull. Acad. Sinica \textbf{11} (1983), 329--336.

\bibitem{Ch} S.\ S.\ Chern, \emph{On the minimal immersions of the two-sphere
in a space of constant curvature.} Problems in Analysis, 27-40. Princeton:
University Press 1970.

\bibitem{DF} M. Dajczer and L. Florit, \emph{A Class of austere submanifolds.}
Illinois Math. J. \textbf{45} (2001), 735--755.

\bibitem{DG} M. Dajczer and D. Gromoll, \emph{The Weierstrass 
representation for complete minimal real Kaehler
submanifolds.} Invent. Math. \textbf{119} (1995), 235--242.

\bibitem{DT} M. Dajczer and R. Tojeiro, \emph{All superconformal surfaces in $\R^4$ in 
terms of minimal surfaces.} Math. Z.  \textbf{261} (2009), 869--890.

\bibitem{DV} M. Dajczer and Th. Vlachos,  \emph{The dual superconformal surface.}\\ 
See http://arxiv.org/abs/1401.1291

\bibitem{Ej} N. Ejiri, \emph{Willmore surfaces with a duality in $\Sf^n(1)$.} 
Proc. London Math. Soc. \textbf{57} (1988), 383--416.

\bibitem{en} ``Enzyklopadie der mathematischen Wissenschaften mit ihrer Anwendungen'' III D 2.1 Teubner.
Leipzig, 1950.

\bibitem{GR} I. Guadalupe and L. Rodr\'iguez, \emph{Normal curvature of surfaces 
in space forms.}\\
Pacific J. Math.  \textbf{106}  (1983), 95--103.

\bibitem{Hu} E. Hulett. \emph{Harmonic superconformal maps of surfaces in $H^n$.} 
J. Geom. Phys. \textbf{42} (2002), 139--165.

\bibitem{Mi} R. Miyaoka, \emph{The family of isometric superconformal
harmonic maps and the affine Toda equations.} J. Reine Angew. Math.
\textbf{481} (1996), 1--25.

\bibitem{mw1} C. Moore and E. Wilson, \emph{A general theory of surfaces.} 
J. Nat. Acad Proc. \textbf{2} (1916), 273--278.

\bibitem{mw2} C. Moore and E. Wilson, \emph{Differential geometry of
two-dimensional surfaces in hyperspaces.} Proc. of the Academy of
Arts and Sciences,  \textbf{52} (1916), 267--368.

\bibitem{Ro} B. Rouxel, \emph{Some geometrical properties of superconformal and
superminimal surfaces in $\mathbb{E}^4$.} Proceedings of the Conference RIGA 2011
Riemannian Geometry and Applications Bucharest, Romania (2011), 255--260.

\bibitem{so} G. Solmon ``A treatise on the analytic geometry of three dimensions''.
New York, Chelsea, 1958-65.

\bibitem{Sp} M. Spivak,``A Comprehensive Introduction to Differential
Geometry''. Vol. IV. Berkeley: Publish or Perish, 1979.

\bibitem{Vl} Th. Vlachos, \emph{Minimal surfaces, Hopf differentials and the
Ricci condition.} Manuscripta Math. \textbf{126} (2008), 201--230.

\bibitem{Wi} P. Wintgen, \emph{Sur l'inegalit\'e de Chen-Willmore.} 
C. R. Acad. Sci. Paris T. Ser. A \textbf{288} (1979), 993--995.
\end{thebibliography}
\end{document}